\documentclass{amsart}
\usepackage{palatino,hhline, bbm}
\usepackage{amsthm, amsmath}
\usepackage{amssymb} 
\usepackage{amscd}
\usepackage{mathrsfs}
\renewcommand{\mathcal}{\mathscr}
\usepackage[hmargin=3cm, vmargin=3.1cm]{geometry}
\usepackage[breaklinks,bookmarksopen,bookmarksnumbered]{hyperref}
\usepackage[all]{xy}

\theoremstyle{plain}
\newtheorem*{thm*}{Theorem}

\theoremstyle{remark}

\newcommand\pr{\noindent\textit{Proof} : }

\newcommand\rond{\kern 1pt{\scriptstyle\circ}\kern 1pt}

\newcommand\Pic{\operatorname{Pic}}

\renewcommand\P{\mathbb{P}}

\renewcommand\O{\mathcal{O}}

\newcommand\iso{\vbox{\hbox to .8cm{\hfill{$\scriptstyle\sim$}\hfill}
\nointerlineskip\hbox to .8cm{{\hfill$\longrightarrow $\hfill}} }}
\newcommand\bir{\vbox{\hbox to .8cm{\hfill{$\scriptstyle\sim$}\hfill}
\nointerlineskip\hbox to .8cm{{\hfill$\dasharrow $\hfill}} }}

\begin{document}
\title{Ulrich bundles on abelian surfaces}
\author[Arnaud Beauville]{Arnaud Beauville}
\address{Laboratoire J.-A. Dieudonn\'e\\
UMR 7351 du CNRS\\
Universit\'e de Nice\\
Parc Valrose\\
F-06108 Nice cedex 2, France}
\email{arnaud.beauville@unice.fr}
 

\begin{abstract}
We prove that any abelian surface admits a rank $2$ Ulrich bundle.
\end{abstract}
\maketitle 
\bigskip
Let $X\subset \P^N$ be a projective variety of dimension $d$ over an algebraically closed field. An \emph{Ulrich bundle} on $X$ is a vector bundle $E$ on $X$ satisfying $H^*(X,E(-1))=\ldots =H^*(X,E(-d))=0$. This notion was introduced in \cite{ES}, where various other characterizations are given; let us just mention that it is equivalent to say that $E$ admits a linear resolution as a $\O_{\P^N}$-module, or that the pushforward of $E$ onto $\P^d$ by a general linear projection is a trivial bundle. 

In \cite{ES} the authors ask whether every projective variety admits an Ulrich bundle. The answer is known only in a few cases: hypersurfaces and complete intersections \cite{HUB}, del Pezzo surfaces \cite[Corollary 6.5]{ES}. 
The case of K3 surfaces  is treated   in \cite{AFO}. In this short note we show that the existence of Ulrich bundles for abelian surfaces follows easily from Serre's construction:
\begin{thm*}
Any abelian surface  $X\subset \P^N$ carries a rank $2$ Ulrich bundle.
\end{thm*}
\pr We put $\dim H^0(X,\O_X(1))=n$.
Let $C$ be a smooth curve in $|\O_X(1)|$; we have $\O_C(1)\cong \omega _C$, and   $g(C)=n+1$. We choose a subset $Z\subset C$ of $n$ general points. Then  $Z$ has the \emph{Cayley-Bacharach property} on $X$ (see for instance \cite{HL},  Theorem 5.1.1): for every $p\in Z$,  any section of $H^0(X,\O_X(1))$  vanishing on $Z\smallsetminus\{p\} $  vanishes on $Z$. Indeed, the image $V$ of the restriction map $H^0(X,\O_X(1))\rightarrow H^0(C,\O_C(1))$ has dimension $n-1$, hence the only element of $V$ vanishing on $n-1$ general points  is zero; thus the only element of $|\O_X(1)|$ containing $Z\smallsetminus\{p\} $ is $C$. 

By \emph{loc. cit.}, there exists a rank 2 vector bundle $E$ on $X$ and an exact sequence
\[0\rightarrow \O_X\xrightarrow{\ s\ } E \longrightarrow \mathcal{I}_Z(1)\rightarrow 0\ .\]
Let $\eta$ be a nonzero element of $\Pic^{\mathrm{o}}(X)$; then $h^0(\omega _C\otimes \eta )=n$, so $H^0(C,\omega _C\otimes \eta(-Z))=0$ since $Z$ is general, and therefore  $H^0(X,\mathcal{I}_Z\eta(1))=0$. Since $\chi (\mathcal{I}_Z\eta(1)))=0$ we have also $H^1(X,\mathcal{I}_Z\eta(1))=0$; from the above exact sequence we conclude that $H^*(X,E\otimes \eta)=0$.

The zero locus of the section $s$ of $E$ is $Z$; since $\det E=\omega _C$,  we get an exact sequence
\[0\rightarrow \O_C(Z)\xrightarrow{\ s_{|C}\ } E_{|C}\longrightarrow \omega _C(-Z)\rightarrow 0\ .\]As above the cohomology of    $\omega _C\otimes \eta(-Z)$ and $\eta(Z)$ vanishes, hence $H^*(C,(E\otimes \eta)_{|C})=0$.
Now from the exact sequence
\[0\rightarrow E(-1)\rightarrow E\rightarrow E_{|C}\rightarrow 0\]we conclude that $H^*(X, E\otimes \eta(-1))=H^*(X,E\otimes \eta)=0$, hence $E\otimes \eta(1)$ is an Ulrich bundle.\qed

\bigskip	
{\small{}\hskip1cm I am indebted to M. Aprodu for pointing out an inaccuracy in the first version of this note.}

\end{document}